\documentclass[preprint,12pt]{elsarticle}

\usepackage{amsmath,amssymb,latexsym,float,epsfig,subfigure}
\usepackage{amsmath, amssymb, amsfonts, amsbsy, latexsym}
\usepackage{amsfonts}
\usepackage{amssymb}
\usepackage{amsthm}
\usepackage{verbatim}
\usepackage{mathrsfs}
\usepackage{amssymb}
 \usepackage{amsthm}

\newcommand{\ud}{\mathrm{d}}
\newcommand{\CHI}{\hbox{\raise .4ex \hbox{$\chi$}}}

\theoremstyle{definition}
\newtheorem{defn}{Definition}
\newtheorem{definition}{Definition}
\newtheorem{lemma}[defn]{Lemma}

\newtheorem{theorem}[defn]{Theorem}

\newtheorem{corollary}[defn]{Corollary}

 \usepackage{lineno}

\journal{Applied and Computational Harmonic Analysis}

\begin{document}

\begin{frontmatter}



\title{Schatten class Fourier Integral Operators}
\author[shan]{Shannon Bishop\corref{cor1}}
\ead{sbishop@math.gatech.edu}
\address[shan]{School of Mathematics, Georgia Institute of Technology, Atlanta, GA 30332 USA}

\begin{abstract}
 Fourier integral operators with sufficiently smooth phase act on the time-frequency content of functions. However time-frequency analysis has only recently been used to analyze these operators.  In this paper, we show that if a Fourier integral operator has a smooth phase function and its symbol is well-localized in time and frequency, then the operator is Schatten \( p \)-class for \( p \in [1,2] \), with inclusion of the symbol in mixed modulation spaces serving as the appropriate measure of time-frequency localization.  Our main results are sharp in the sense that larger mixed modulation spaces necessarily contain symbols of Fourier integral operators that are not Schatten \( p \)-class.

\end{abstract}


\end{frontmatter}


\section{Introduction}

	  Classical Fourier integral operators, which arise in the study of hyperbolic differential equations (see \cite{treves2}), are operators of the form \begin{equation}
Af(x) = \int a(x,\xi) \widehat{f}(\xi) e^{2 \pi i \varphi(x,\xi)} \, \ud \xi \label{eqn:easyfio}.
\end{equation}  In this case \( a  \) is the symbol and \( \Phi \) is the phase function of the operator.  Fourier integral operators generalize pseudodifferential operators in the sense that the pseudodifferential operator with Kohn-Nirenberg symbol \( \sigma \) is the Fourier integral operator with symbol \( \sigma \) and phase function \( \varphi(x, \xi) = x \cdot \xi \).
More generally, an operator of the form 
 \begin{equation}
Af(x) = \iint b(x,y, \xi) f(y) e^{2 \pi i \psi(x,y,\xi)} \, \ud y \, \ud \xi \label{eqn:hardfio}
\end{equation} is also called a Fourier integral operator and \( b \), \( \psi \) are the symbol and phase function of \( A \) respectively.
The properties of Fourier integral operators with smooth symbols and phase functions have been studied extensively.  In particular the boundedness properties of such operators are well-known.   If the symbol and phase function belong to \( C^{\infty} \) and the symbol belongs to an appropriate Hormander symbol class, then the operator is bounded on \( L^{p} \) (see \cite{ruz}, \cite{steinFIO} and the references therein), bounded from \( H^{1} \) to \( L^{1} \) (see \cite{steinFIO} and the generalization in \cite{tao}), bounded on weighted Sobolev spaces (see \cite{ruz2}) and bounded on \( \mathcal{F} L^{p} \) (see \cite{rod1}). More recently, in \cite{candesFIO} and \cite{labateFIO}, it was shown that the curvelet and shearlet representations of a Fourier integral operator of the form (\ref{eqn:easyfio}) are sparse, provided  \( a, \varphi \in C^{\infty} \) and \(a \) is in an appropriate Hormander symbol class.  Much less is known about Fourier integral operators with non-smooth symbols.

However, for one particular type of Fourier integral operator, the pseudodifferential operator, much can be said, even in the case of a non-smooth symbols.  Specifically, time-frequency analysis can be used to describe the properties of pseudodifferential operators, with the modulation spaces (which are Banach spaces characterized by time-frequency shifts and mixed norms) serving as appropriate symbol spaces for studying continuity and Schatten class properties of pseudodifferential operators.  Using Gabor frames, elements in the modulation spaces can be decomposed into a superposition of time-frequency shifts, and this Gabor frame decomposition of the symbol of a pseudodifferential operator can be used to characterize the properties of the operator.  Results of this type appear in \cite{czaja}, \cite{modbounded}, \cite{heilpdo}, \cite{labate}, \cite{tach2} and \cite{toftpdo}.  
	  
	  Both pseudodifferential operators and Fourier integral operators with smooth phase functions act on the time-frequency content of functions, although the time-frequency action of a Fourier integral operator is much more general and less explicit than the action of a pseudodifferential operator.   However, this action still suggests that time-frequency analysis may play an important role in understanding Fourier integral operators with non-smooth symbols.  Indeed, recent results confirm this intuition.  In \cite{boulk} it was shown that inclusion of the symbol of a Fourier integral operator with smooth phase in Sj\"ostrand's class (a space containing non-smooth symbols) implies boundedness of the operator on \( L^{2}(\mathbb{R}^{d}) \).   In \cite{cordero2} and \cite{cordero1}, the authors use time-frequency analysis to prove the boundedness of Fourier integral operators of the form (\ref{eqn:easyfio}) on the time-frequency localization space \( M^{p} \), provided \( a \in C^{2N} \), a much less strenuous condition than the traditional smoothness assumptions on the symbol.   More generally, in \cite{toftFIO1} and \cite{toftFIO2}, the authors prove Schatten \( p \)-class membership for Fourier integral operators with sufficiently smooth phase functions whose symbols belong to \( M^{p,1} \), another space containing non-smooth symbols.   Note that while Fourier integral operators generalize  pseudodifferential operators, pseudodifferential operator analysis techniques do not appear to generalize to Fourier integral operators. The results in \cite{boulk}, \cite{toftFIO1}, \cite{toftFIO2}, \cite{cordero2} and \cite{cordero1} are proved with new Gabor frame techniques.  
	  
	 	  In this paper, we use time-frequency analysis techniques to characterize Schatten class Fourier integral operators, with the mixed modulation spaces, which are natural generalizations of the traditional modulation spaces, serving as appropriate symbol classes.  For operators of the form (\ref{eqn:easyfio}), we show that inclusion of the symbol in \( M(c)^{2,\dots,2, p,\dots,p} \) for certain permutations \( c \)  and \( p \in [1,2] \) almost guarantees the operator is Schatten \(p \)-class when the phase is smooth.  We also prove that if the symbol of a Fourier integral operator of the form (\ref{eqn:hardfio}) belongs to the mixed modulation space \( M(c)^{2, \dots, 2,p, \dots, p,1,\infty} \) or \( M(c)^{\infty, 2, \dots, 2,p, \dots, p,1} \) for appropriate permutations \( c \) and if the phase function is sufficiently smooth, then the operator is Schatten \( p \)-class for \( p \in [1,2] \).  Our results for Fourier integral operators of the form (\ref{eqn:easyfio}) improve upon existing results of this type and are sharp in the sense that mixed modulation spaces larger than \( M(c)^{2,\dots,2, p,\dots,p} \) contain Fourier integral operators that are not Schatten \( p \)-class.  Although our results for Fourier integral operators of the form (\ref{eqn:hardfio}) are not directly comparable to previously known Schatten class results for Fourier integral operators, namely those in \cite{toftFIO1} and \cite{toftFIO2}, they seem stronger in the sense that \( M(c)^{2, \dots, 2,p, \dots, p,1,\infty} \), \( M(c)^{\infty, 2, \dots, 2,p, \dots, p,1} \) are isomorphic to \( \ell^{2, \dots, 2,p, \dots, p,1,\infty} \) and \( \ell^{\infty, 2, \dots, 2,p, \dots, p,1} \), respectively,  while  \( M^{p,1} \) is isomorphic to \( \ell^{p, \dots, p,1,\dots, 1} \) and \( \ell^{p, \dots, p,1,\dots, 1} \subsetneq \ell^{2, \dots, 2,p, \dots, p,1,\infty} \) and \( \ell^{p, \dots, p,1,\dots, 1} \subsetneq \ell^{\infty, 2, \dots, 2,p, \dots, p,1} \).  Furthermore, our main results  for Fourier integral operators of the form (\ref{eqn:hardfio}) are sharp in the sense that larger mixed modulation spaces  contain symbols of Fourier integral operators that are not Schatten \( p \)-class.   However, the results in this paper are not applicable to the special case of Fourier integral operators of the form
	 	   \begin{equation}
Af(x) = \iint b(x, \xi) f(y) e^{2 \pi i \psi(x, y, \xi)} \, \ud y \, \ud \xi,
\end{equation}
whereas \cite{toftFIO1} and \cite{toftFIO2} do give Schatten class results for operators of this type.

The remainder of this paper is organized as follows.  Section 2 contains background information.  In Section 3, we prove a time-frequency condition on the product of the symbol and phase of a Fourier integral operator that ensures the operator is Schatten class.  In Section 4, we use mixed modulation space product embeddings to find sharp time-frequency conditions on the symbol of a Fourier integral operator that ensure the operator is Schatten class.  

\section{Background}

\subsection{Mixed norm spaces}

\begin{definition}

 Given measure spaces \( \left( X_{i}, \mu_{i} \right) \) and indices \( p_{i} \in [1, \infty] \) for \( i = 1, 2, \dots, d \),  we let  \[ L^{p_{1}, p_{2}, \dots, p_{d}} \left( X_{1}, X_{2}, \dots, X_{d}, \mu_{1}, \mu_{2}, \dots,  \mu_{d}\right) \] consist of all measurable functions \( F: X_{1} \times X_{2} \times \cdots \times X_{d} \to \mathbb{C} \) for which the following norm is finite:
\[ \left\| F \right\|_{L^{p_{1}, p_{2}, \dots, p_{d}}}  = \biggl( \int_{X_{d}} \dots   \biggl( \int_{X_{1}} \left| F \left(x_{1}, \dots, x_{d} \right)   \right|^{p_{1}}   \, \ud \mu_{1} (x_{1})  \biggr)^{\frac{p_{2}}{p_{1}}} \cdots  \, \ud \mu_{d} (x_{d}) \biggr)^{\frac{1}{p_{d}}}, \]
with the usual modifications for those indices \( p_{i} \) which equal \( \infty \).  

  If \( X_{i} = \mathbb{R} \) and \(  \mu_{i} \) is Lebesgue measure on \( \mathbb{R} \) for all \( i = 1, 2, \dots, d \), then we simply write \( L^{p_{1}, p_{2}, \dots, p_{d}} \).  If each \( X_{i} \) is countable and \( \mu_{i} \) is counting measure on \( X_{i} \) we simply write \( \ell^{p_{1}, p_{2}, \dots, p_{d}}_{w} \left( X_{1}, X_{2}, \dots, X_{d}\right) \).
\end{definition}

The mixed norm spaces \( L^{p_{1}, p_{2}, \dots, p_{d}} \left( X_{1}, X_{2}, \dots, X_{d}, \mu_{1}, \mu_{2}, \dots,  \mu_{2d}\right) \) are generalizations of the classical spaces \( L^{p} \), and the proof that \( L^{p} \) is a Banach space can be extended to the mixed norm spaces (see \cite{oldie}).

A Wiener amalgam norm is a type of mixed norm that measures local boundedness with global decay.
\begin{definition}  Suppose \( p_{1}, \dots, p_{d} \in [1, \infty] \).  Define a norm by
\[   \left\| f \right\|_{W ( L^{1}(\mathbb{R}^{d}))} = \sum_{n \in \mathbb{Z}^{d}} \left\| f \CHI_{[0,1]^{d} + n } \right\|_{\infty}.\]  The \emph{Wiener amalgam space} \(W ( L^{1}(\mathbb{R}^{d})) \) is the set of functions for which this norm is finite.   
\end{definition}

For any multi-index \( \alpha = \left( \alpha_{1}, \dots, \alpha_{d} \right) \in \mathbb{R}^{d} \) with \(  \alpha_{1}, \dots, \alpha_{d} \in (0, \infty) \), we can define an equivalent norm on \( W ( L^{1}(\mathbb{R}^{d})) \) by \begin{equation} \label{eqn:star}   \sum_{n \in \mathbb{Z}^{d}} \left\| f \CHI_{ \alpha \cdot [0,1]^{d} + \alpha \cdot n } \right\|_{\infty}. \end{equation}

\subsection{Gabor Transform}

Suppose \( f:\mathbb{R}^{d} \to \mathbb{C} \) is measurable.  For \( x, \xi \in \mathbb{R}^{d} \) define the translation operator \( T_{x} \) and modulation operator \( M_{\xi} \) by
\[ T_{x}f(t) = f( t -x ) \hspace{2pc} \textrm{ and } \hspace{2pc}  M_{\xi}f(t) = e^{2 \pi i t \cdot \xi} f(t),\] and define the time-frequency shift \( \pi_{(x, \xi)} \) by \( \pi_{(x, \xi)} = M_{\xi} T_{x} \).

\begin{definition}  Fix \( \phi \in \mathscr{S}(\mathbb{R}^{d}) \).  Given \( f \in \mathscr{S}'(\mathbb{R}^{d}) \), the \emph{Gabor transform} or \emph{short-time Fourier transform} of \( f \) with respect to \( \phi \) is 
\[ V_{\phi} f( x, \xi ) = \int_{\mathbb{R}^{d}} f(t) \overline{ \phi(t-x)} e^{-2 \pi i \xi \cdot t} \, \ud t = \langle f, M_{\xi} T_{ x} \phi \rangle,  \hspace{2pc} x, \xi \in \mathbb{R}^{d}. \]  The function \( \phi \) is called the \emph{window function} of the Gabor transform.
\end{definition}

The value of \( V_{\phi} f( x, \xi ) \) gives information about the time-frequency content of \( f \) around \( x \) in time and \( \xi \) in frequency.  See \cite{groch} for background and information about the Gabor transform.

If \( c \) is a permutation of \( \left\{ 1,2,\dots, 2d\right\} \), we identify  \(c \)  with the bijection \( \mathfrak{c}: \mathbb{R}^{2d} \to \mathbb{R}^{2d} \) given by \( \mathfrak{c}(x_{1},\dots x_{2d}) = (x_{c(1)}, \dots, x_{c(2d)}) \).

\begin{definition}  Suppose \( \phi \in \mathscr{S}(\mathbb{R}^{d}) \) and \( c \) is a permutation of \( \left\{ 1,2,\dots, 2d\right\}  \) corresponding to the map \( \mathfrak{c} \).  Let \( M(c)^{p_{1}, p_{2}, \dots, p_{2d}} \) be the \emph{mixed modulation space} consisting of all \(  f \in \mathscr{S}'(\mathbb{R}^{d}) \) for which
\[ \left\| f \right\|_{M(c)^{p_{1}, p_{2}, \dots, p_{2d}}} = \left\| V_{\phi} f \circ \mathfrak{c} \right\|_{L^{p_{1}, p_{2}, \dots, p_{2d}}} < \infty. \] 
\end{definition}

The mixed modulation spaces generalize the traditional modulation spaces \( M^{p,q}(\mathbb{R}^{d}) \).  In particular, if \( c \) is the identity permutation  and \( p = p_{1} = p_{2} = \dots = p_{d} \) and \( q = p_{d+1} = \dots = p_{2d} \) then \( M(c)^{p_{1}, p_{2}, \dots, p_{2d}} = M^{p,q}(\mathbb{R}^{d})  \).  Furthermore, the most interesting properties of modulation spaces carry over to the mixed modulation spaces.  See \cite{mine} for more information about mixed modulation spaces. 

\subsection{Gabor Frames}

\begin{definition} Suppose \( \Lambda \subset \mathbb{R}^{2d} \) is a countable set.  A \emph{Gabor frame} for a Hilbert space \( L^{2}(\mathbb{R}^{d}) \) is a set \( \left\{ M_{\xi} T_{x} \phi \right\}_{(x, \xi) \in \Lambda }  \subset L^{2}(\mathbb{R}^{d}) \) such that there are \( A, B > 0 \) with 
\[ A \left\| f \right\|^{2} \leq \sum_{(x, \xi) \in \Lambda} \left| \langle f, M_{\xi} T_{x} \rangle \right|^{2} \leq B \left\| f \right\|^{2}\] for all \( f \in L^{2}(\mathbb{R}^{d}) \).  In this case \( A, B \) are \emph{frame bounds}.  If we can take \( A = B \) then 
\( \left\{ M_{\xi} T_{x} \phi \right\}_{(x, \xi) \in \Lambda} \)  is a \emph{tight frame}.  A tight frame is \emph{Parseval} if we can choose \( A = B = 1 \).  
\end{definition}

Gabor frames give nonorthogonal expansions of elements of \( L^{2}(\mathbb{R}^{d}) \) in terms of the frame elements, and  these expansions are  stable but usually redundant.  If  \( \left\{ M_{\xi} T_{x} \phi \right\}_{(x, \xi) \in \Lambda} \) is a frame for \( L^{2}(\mathbb{R}^{d}) \), there is a dual sequence \( \{ \tilde{\phi}_{x, \xi} \}_{(x, \xi) \in \Lambda} \subset L^{2}(\mathbb{R}^{d}) \) such that 
\[ f = \sum_{(x, \xi) \in \Lambda}  \langle f, M_{\xi} T_{x} \phi \rangle  \tilde{\phi}_{x, \xi}  = \sum_{(x,\xi) \in \Lambda}  \langle f, \tilde{\phi}_{x, \xi} \rangle M_{\xi} T_{x} \phi  \] 
for all \( f \in L^{2}(\mathbb{R}^{d}) \), and the sequence \( \{ \tilde{\phi}_{x, \xi} \}_{(x, \xi) \in \Lambda} \) can be chosen to be a frame for \( L^{2}(\mathbb{R}^{d})  \). In particular, if \( \left\{M_{\xi} T_{x} \phi \right\}_{x \in X} \) is a tight frame for \(  L^{2}(\mathbb{R}^{d}) \) with frame bound \( B \), we have
\[ f = B^{-1} \sum_{(x, \xi) \in \Lambda} \langle f , M_{\xi} T_{x} \phi \rangle M_{\xi} T_{x} \phi \hspace{2pc} \forall f \in L^{2}(\mathbb{R}^{d}).\] 
See \cite{lsu} for general background on Gabor frames. 

A Gabor frame with a nice generator and regular index set has powerful properties in time-frequency spaces beyond \( L^{2}(\mathbb{R}^{d}) \).  The following theorem from \cite{mine} makes this idea precise.

\begin{theorem}  \label{theorem:banachframe} Fix \( \beta > 0 \).  Suppose \(  p_{1}, p_{2}, \dots, p_{2d} \in [1, \infty] \) and \( \psi \in M(c)^{1,\dots,1} \).  Further suppose that \( \left\{ \pi_{\beta n} \psi \right\}_{n \in \mathbb{Z}^{2d}} \) is a frame for \( L^{2}(\mathbb{R}^{d}) \) with dual frame \( \left\{ \pi_{\beta n} \gamma \right\}_{n \in \mathbb{Z}^{2d}} \). Then
\begin{itemize}
\item[(a)] \( \left\{ \pi_{\beta n} \psi \right\}_{n \in \mathbb{Z}^{2d}} \) is a Banach frame for \(  M(c)^{p_{1}, p_{2}, \dots, p_{2d}} \) and there exist \(0 < A \leq B < \infty \) independent of \(   p_{1}, p_{2}, \dots, p_{2d} \) such that \[ A \left\| f \right\|_{ M(c)^{p_{1}, p_{2}, \dots, p_{2d}}} \leq \left\| V_{\psi} f \circ \mathfrak{c} \big|_{\beta \mathbb{Z}^{2d}} \right\|_{\ell^{p_{1}, p_{2}, \dots, p_{2d}}} \leq B \left\| f \right\|_{ M(c)^{p_{1}, p_{2}, \dots, p_{2d}}}, \] for all \( f \in M(c)^{p_{1}, p_{2}, \dots, p_{2d}} \).
\item[(b)]  If \( p_{1}, p_{2}, \dots, p_{2d} \in [1, \infty) \) then \begin{displaymath} f = \sum_{m \in \mathbb{Z}^{2d}} \langle f, \pi_{\beta m} \psi \rangle \, \pi_{\beta m} \gamma =  \sum_{m \in \mathbb{Z}^{2d}} \langle f, \pi_{\beta m} \gamma \rangle \, \pi_{\beta m} \psi \end{displaymath} for all \( f \in M(c)^{p_{1}, p_{2}, \dots, p_{2d}} \) with unconditional convergence in \( M(c)^{p_{1}, p_{2}, \dots, p_{2d}} \). \smallskip
\item[(c)]  If \( p_{1}, p_{2}, \dots, p_{2d} \in [1, \infty] \) then  \begin{displaymath} f = \sum_{m \in \mathbb{Z}^{2d}} \langle f, \pi_{\beta m} \psi \rangle \, \pi_{\beta m} \gamma =  \sum_{m \in \mathbb{Z}^{2d}} \langle f, \pi_{\beta m} \gamma \rangle \, \pi_{\beta m} \psi \end{displaymath} for all \( f \in M(c)^{p_{1}, p_{2}, \dots, p_{2d}} \) with weak* convergence in \( M(c)^{\infty, \dots, \infty} \).
\end{itemize}
\end{theorem}

 \subsection{Permutations}
 In Section 2 and 3 certain mixed modulation spaces will be useful in classifying Schatten class Fourier integral operators.  These mixed modulation spaces depend on the permutations of the variables of the Gabor transform defined as follows.

 \begin{definition}    Suppose \( c \) is a permutation on \( \left\{ 1,2,\dots, 4d \right\} \).
A \emph{first slice} permutation \( c \) is one that satisfies
\begin{itemize}
\item[(a)] \( c \) maps \( \left\{ 1,2,\dots, d, 2d+1, 2d+2, \dots, 3d \right\}\) to  \( \left\{ 1,2,\dots, 2d \right\} \)  bijectively and
\item[(b)] \( c \) maps  \( \left\{d+1,d+2, \dots, 2d, 3d+1, 3d+2, \dots, 4d\right\} \) to \( \left\{ 2d+1, \dots, 4d \right\} \)  bijectively.
\end{itemize}
A \emph{second slice} permutation \( c \) is one that satisfies
\begin{itemize}
\item[(a)] \( c \) maps \( \left\{ d+1,d+2,\dots, 2d, 3d+1, 3d+2, \dots, 4d \right\}\) to  \( \left\{ 1,2,\dots, 2d \right\} \)  bijectively and
\item[(b)] \( c \) maps  \( \left\{1,2, \dots, d, 2d+1, 2d+2, \dots, 3d\right\} \) to \( \left\{ 2d+1, \dots, 4d \right\} \)  bijectively.
\end{itemize}
\end{definition}

\begin{definition}
Suppose \( c \) is a permutation on \( \left\{ 1,2,\dots, 6d \right\} \).
A \emph{first FIO slice permutation} \( c \) is a permutation  of \( \left\{ 1,2,\dots, 6d \right\} \) such that
\begin{itemize}
\item[(a)] \( c \) maps \( \left\{ 1, 2, \dots, d, 3d+1, 3d+2, \dots, 4d \right\} \) to \( \left\{ 1, 2, \dots, 2d \right\} \),
\item[(b)] \( c \) maps \( \left\{ d+1, d+2, \dots, 2d, 4d+1, 4d+2,\dots, 5d \right\} \) to \( \left\{ 2d+1, 2d+2, \dots, 4d \right\} \),
\item[(c)] \( c \) maps \( \left\{ 2d+1, 2d+2,\dots, 3d \right\} \) to \( \left\{4d+1, 4d+2,\dots, 5d \right\} \), and
\item[(d)] \( c \) maps \( \left\{ 5d+1,5d+2,\dots,6d \right\} \) to \( \left\{5d+1, 5d+2,\dots, 6d \right\} \).
\end{itemize}
A \emph{second FIO slice permutation} \( c \) is a permutation  of \( \left\{ 1,2,\dots, 6d \right\} \) such that
\begin{itemize}
\item[(a)] \( c \) maps \( \left\{ d+1, d+2, \dots, 2d, 4d+1, 4d+2, \dots, 5d \right\} \) to \( \left\{ 1, 2, \dots, 2d \right\} \),
\item[(b)] \( c \) maps \( \left\{ 1, 2, \dots, d, 3d+1, 3d+2,\dots, 4d \right\} \) to \( \left\{ 2d+1, 2d+2, \dots, 4d \right\} \),
\item[(c)] \( c \) maps \( \left\{ 2d+1, 2d+2,\dots, 3d \right\} \) to \( \left\{4d+1, 4d+2,\dots, 5d \right\} \), and
\item[(d)] \( c \) maps \( \left\{ 5d+1,5d+2,\dots,6d \right\} \) to \( \left\{5d+1, 5d+2,\dots, 6d \right\} \).
\end{itemize}
A \emph{first FIO symbol permutation} \( c \) is a permutation  of \( \left\{ 1,2,\dots, 6d \right\} \) such that
\begin{itemize}
\item[(a)] \( c \) maps \( \left\{ 5d+1, 5d+2, \dots, 6d \right\} \) to \( \left\{ 1, 2, \dots, d \right\} \),
\item[(b)] \( c \) maps \( \left\{ 1, 2, \dots, d, 3d+1, 3d+2,\dots, 4d \right\} \) to \( \left\{ d+1, d+2, \dots, 3d \right\} \),
\item[(c)] \( c \) maps \( \left\{ d+1,\dots, 2d, 4d+1, 4d+2, \dots, 5d \right\} \) to \( \left\{3d+1, 3d+2,\dots, 5d \right\} \), and
\item[(d)] \( c \) maps \( \left\{ 2d+1,2d+2,\dots,3d \right\} \) to \( \left\{5d+1, 5d+2,\dots, 6d \right\} \).
\end{itemize}
A \emph{second FIO symbol permutation} \( c \) is a permutation  of \( \left\{ 1,2,\dots, 6d \right\} \) such that
\begin{itemize}
\item[(a)] \( c \) maps \( \left\{ 5d+1, 5d+2, \dots, 6d \right\} \) to \( \left\{ 1, 2, \dots, d \right\} \),
\item[(b)] \( c \) maps  \( \left\{ d+1, d+2,\dots, 2d, 4d+1, 4d+2, \dots, 5d \right\} \)  to \( \left\{ d+1, d+2, \dots, 3d \right\} \),
\item[(c)] \( c \) maps \( \left\{ 1, 2, \dots, d, 3d+1, 3d+2,\dots, 4d \right\} \) to \( \left\{3d+1, 3d+2,\dots, 5d \right\} \), and
\item[(d)] \( c \) maps \( \left\{ 2d+1,2d+2,\dots,3d \right\} \) to \( \left\{5d+1, 5d+2,\dots, 6d \right\} \).
\end{itemize}
\end{definition}

 \subsection{Schatten class Operators}  
 \begin{definition} Fix \(1\leq p  < \infty\).  Suppose \( A:  L^{2}(\mathbb{R}^{d})  \to  L^{2}(\mathbb{R}^{d})  \) is a linear operator.  We say \( A \) is \emph{Schatten \( p \)-class} and write \( A \in \mathcal{I}_{p}( L^{2}(\mathbb{R}^{d}) ) \) if \[ \left\| A \right\|_{\mathcal{I}_{p}} = \sup \biggl( \sum_{n \in \mathbb{N}} \left| \langle A f_{n}, g_{n} \rangle \right|^{p} \biggr)^{\frac{1}{p}} < \infty,\] where the supremum is taken over all pairs of orthonormal sequences \( \left\{ f_{n} \right\}_{n \in \mathbb{N}} \),  \( \left\{ g_{n} \right\}_{n \in \mathbb{N}} \) in \(  L^{2}(\mathbb{R}^{d})  \).
\end{definition}

Equivalently, an operator is Schatten \( p\)-class if its singular values constitute an \( \ell^{p} \) sequence.  Consequently, trace-class operators are exactly the Schatten 1-class operators and Hilbert-Schmidt operators are the Schatten 2-class operators.  For \(p = \infty \), we define Schatten \( p\)-class operators to be bounded operators.

\subsection{Integral Operators}

An integral operator with kernel \(k \) is one of the form 
\[ Af(x) = \int k(x,y) f(y) \, \ud y.\] 

If the kernel of an integral operator has sufficient time-frequency concentration, then the operator is Schatten class.  This is made precise in the following theorem from \cite{mine}.

\begin{theorem} \label{thm:intop} Suppose \( c \) is a first or second slice permutation on \( \left\{1,2, \cdots, 4d \right\} \) and \( p_{1} = p_{2} = \cdots = p_{2d} = 2 \) and \( p_{2d+1} = \cdots = p_{4d} = p \) for some \( p \in [1,2] \).  If \( A \) is an integral operator with kernel \( k \) and \( k \in M(c)^{p_{1}, p_{2}, \cdots, p_{4d}} \), then \( A \in \mathcal{I}_{p}( L^{2}(\mathbb{R}^{d}) ) \).
Furthermore, this result is sharp in the sense that if one of the following conditions holds then there are integral operators not in \( \mathcal{I}_{p}(L^{2}(\mathbb{R}^{d})) \) whose kernel lies in \( M(c)^{q_{1},q_{2},\dots, q_{4d}} \).
\begin{itemize}
\item[(a)]  At least one of \(  q_{1} ,  \dots  , q_{2d} \) is larger than 2.
\item[(b)] At least one of \(  q_{2d+1} , \dots , q_{4d} \) is larger than \( p \).
\end{itemize} 
\end{theorem}

\section{Schatten class results for Fourier Integral Operators}

Notice that if \( A \) is the Fourier integral operator in (\ref{eqn:easyfio}), then \( A = T \mathcal{F} \), where \( \mathcal{F} \) is the Fourier transform and \( T \) is the integral operator with kernel \( k (x,y) = a(x, y)  e^{2 \pi i \varphi(x,y)} \).  Hence, we obtain the following theorem as a direct consequence of Theorem \ref{thm:intop}.

\begin{theorem}  \label{thm:FIOlemma} Suppose \( A \) is Fourier integral operator of the form (\ref{eqn:easyfio}). Let \( c \) be a first or second slice permutation on \( \left\{1,2, \cdots, 4d \right\} \) and \( p_{1} = p_{2} = \cdots = p_{2d} = 2 \) and \( p_{2d+1} = \cdots = p_{4d} = p \) for some \( p \in [1,2] \). If \( a e^{2 \pi i \varphi} \in M(c)^{p_{1}, p_{2}, \cdots, p_{4d}} \), then \( A \in \mathcal{I}_{p}( L^{2}(\mathbb{R}^{d}) ) \).
\end{theorem}

Although Theorem \ref{thm:intop} is not immediately applicable to the more general Fourier integral operators of the form (\ref{eqn:hardfio}), it can be adapted to this type of operator to give time-frequency conditions on the product of the symbol and phase which ensure the operator is Schatten class, as in the following theorem.

\begin{theorem} \label{thm:FIOmainone}  Suppose \( A \) is a Fourier integral operator of the form (\ref{eqn:hardfio})  with symbol \( b \) and phase function \( \psi \).  Suppose \( p \in [1,2] \) and \( c \) is a first or second FIO slice permutation. Let \( p_{1} = p_{2} = \dots = p_{2d} = 2 \), \( p_{2d+1} = p_{2d+2} = \dots = p_{4d} = p \), \( p_{4d+1} = p_{4d+2} = \dots = p_{5d}= 1 \) and \( p_{5d+1} = p_{5d+2} = \dots = p_{6d} = \infty \). If \( b e^{2 \pi i \psi} \in M(c)^{p_{1},p_{2},\dots,p_{6d}} \), then \( A \in \mathcal{I}_{p}(L^{2}(\mathbb{R}^{d})) \).
\end{theorem}
\begin{proof}  We prove the result in the case \( c \) is a second FIO slice permutation.  The case that \( c \) is a first FIO slice permutation can be proven similarly.  

Let \( \left\{ f_{n} \right\}_{n \in \mathbb{N}} \), \( \left\{ g_{n} \right\}_{n \in \mathbb{N}} \) be arbitrary orthonormal sequences in \( L^{2}(\mathbb{R}^{d}) \) and let  \( \left\{ M_{\alpha k_{2}} T_{\alpha k_{1}} \phi \right\}_{k_{1},k_{2} \in \mathbb{Z}^{d}} \) be a Parseval Gabor frame for \( L^{2}(\mathbb{R}^{d}) \) with \( \phi \in M^{1,1}(\mathbb{R}^{d}) \) and \( \alpha > 0 \). By Proposition 12.1.4 in \cite{groch}, \( \phi \in M^{1,1}(\mathbb{R}^{d}) \) implies \( \hat{\phi} \in W( L^{1}(\mathbb{R}^{d}))\). Using the equivalent norm given in (\ref{eqn:star}), there is a  \( C \) such that 
\[ \left\| \left\{ \left\| \hat{\phi} \CHI_{ \alpha  [0,1]^{d} + \alpha  n } \right\|_{\infty} \right\}_{n \in \mathbb{Z}^{d}} \right\|_{\ell^{1}(\mathbb{Z}^{d})} \leq C \left\| \hat{\phi} \right\|_{W ( L^{1}(\mathbb{R}^{d}) )}.\]

 Using the definition of \( A \) we have 
\[ \langle A f_{n}, g_{n} \rangle = \langle b e^{2 \pi i \psi}, g_{n} \otimes \overline{f_{n}} \otimes 1 \rangle. \]   Since \( 1 \in M^{\infty, 1}(\mathbb{R}^{d}) \), Theorem \ref{theorem:banachframe} implies \[ 1 = \sum_{k_{1}, k_{2} \in \mathbb{Z}^{d}} \langle 1,  M_{\alpha k_{2}} T_{\alpha k_{1}} \phi \rangle  M_{\alpha k_{2}} T_{\alpha k_{1}} \phi \hspace{1pc} \textrm{ weakly.}  \] Thus
\begin{align*}
 \langle A f_{n}, g_{n} \rangle & = \sum_{k_{1},k_{2} \in \mathbb{Z}^{d}} \overline{\langle 1, M_{\alpha k_{2}} T_{\alpha k_{1}} \phi \rangle} \langle b e^{2 \pi i \psi}, g_{n}\otimes \overline{f_{n}} \otimes  M_{\alpha k_{2}} T_{\alpha k_{1}} \phi \rangle \\
& =  \sum_{k_{1},k_{2}\in \mathbb{Z}^{d}} \overline{\langle 1, M_{\alpha k_{2}} T_{\alpha k_{1}} \phi \rangle} \langle A_{k_{1},k_{2}} f_{n}, g_{n} \rangle,
\end{align*}
where \( A_{k_{1},k_{2}} \) is the integral operator with kernel
\[ k_{k_{1},k_{2}}(x,y) = \int b(x,y,\xi) e^{2 \pi i \psi(x,y,\xi)} \overline{ M_{\alpha k_{2}} T_{\alpha k_{1}} \phi(\xi)}  \, \ud \xi. \]

In the case \( p = 1 \) we have
\begin{align*}
\sum_{n \in \mathbb{N}} \left| \langle A f_{n}, g_{n} \rangle \right| & = \sum_{n \in \mathbb{N}} \left| \sum_{k_{1},k_{2} \in \mathbb{Z}^{d}} \overline{\langle 1, M_{\alpha k_{2}} T_{\alpha k_{1}} \phi \rangle} \langle A_{k_{1},k_{2}} f_{n}, g_{n} \rangle \right| \\
& \leq  \sum_{k_{1},k_{2} \in \mathbb{Z}^{d}} \sum_{n \in \mathbb{N}}  \left| \overline{\langle 1, M_{\alpha k_{2}} T_{\alpha k_{1}} \phi \rangle} \langle A_{k_{1},k_{2}} f_{n}, g_{n} \rangle \right| \\
& =  \sum_{k_{1},k_{2} \in \mathbb{Z}^{d}} \sum_{n \in \mathbb{N}}  \left| \hat{\phi}(\alpha k_{2}) \right| \left| \langle A_{k_{1},k_{2}} f_{n}, g_{n} \rangle \right| \\
& = \sum_{k_{1},k_{2} \in \mathbb{Z}^{d}}  \left| \hat{\phi}(\alpha k_{2}) \right|\sum_{n \in \mathbb{N}}  \left| \langle A_{k_{1},k_{2}} f_{n}, g_{n} \rangle \right| \\
& \leq \left( \sum_{k_{2} \in \mathbb{Z}^{d}} \left| \hat{\phi}(\alpha k_{2}) \right| \right) \left( \sup_{k_{2} \in \mathbb{Z}^{d}} \sum_{k_{1} \in \mathbb{Z}^{d}}\sum_{n \in \mathbb{N}}  \left| \langle A_{k_{1},k_{2}} f_{n}, g_{n} \rangle \right| \right) \\
& \leq C \left\| \hat{\phi} \right\|_{W(L^{1}(\mathbb{R}^{d}))} \sup_{k_{2} \in \mathbb{Z}^{d}} \sum_{k_{1} \in \mathbb{Z}^{d}} \sum_{n \in \mathbb{N}}  \left| \langle A_{k_{1},k_{2}} f_{n}, g_{n} \rangle \right| \\
& \leq C  \left\| \hat{\phi} \right\|_{W(L^{1}(\mathbb{R}^{d}))} \sup_{k_{2} \in \mathbb{Z}^{d}} \sum_{k_{1} \in \mathbb{Z}^{d}} \left\| A_{k_{1},k_{2}} \right\|_{\mathcal{I}_{1}}
\end{align*}
By Theorem \ref{thm:intop}, we have
\begin{align*}
\left\|  A_{k_{1},k_{2}} \right\|_{\mathcal{I}_{1}} & \leq \sum_{n_{1},n_{2} \in \mathbb{Z}^{d}} \left( \sum_{m_{1}, m_{2} \in \mathbb{Z}^{d}} \left| V_{\phi \otimes \phi} k_{k_{1},k_{2}} (  \alpha n_{1},\alpha m_{1}, \alpha n_{2},\alpha m_{2}) \right|^{2} \right)^{\frac{1}{2}} \\
& =  \sum_{n_{1},n_{2}\in \mathbb{Z}^{d}} \left( \sum_{m_{1}, m_{2}\in \mathbb{Z}^{d}} \left| V_{\Phi} (b e^{2 \pi i \psi})( \alpha n_{1},\alpha m_{1}, \alpha k_{1}, \alpha n_{2},\alpha m_{2}, \alpha k_{2}) \right|^{2} \right)^{\frac{1}{2}},
\end{align*}
where \( \Phi = \phi \otimes \phi \otimes \phi \). Thus if
\[ \sup_{k_{2} \in \mathbb{Z}^{d}} \sum_{k_{1} \in \mathbb{Z}^{d}} \sum_{n_{1},n_{2} \in \mathbb{Z}^{d}} \left( \sum_{m_{1}, m_{2} \in \mathbb{Z}^{d}} \left| V_{\Phi} (b e^{2 \pi i \psi})(\alpha n_{1},\alpha m_{1}, \alpha k_{1}, \alpha n_{2},\alpha m_{2}, \alpha k_{2}) \right|^{2} \right)^{\frac{1}{2}} < \infty, \] then  \( A \in \mathcal{I}_{1}(L^{2}(\mathbb{R}^{d})) \).  Notice that this quantity is finite if and only if \( b e^{2 \pi i \psi} \in M(c)^{p_{1},p_{2},\dots,p_{6d}} \).

For the case \( p = 2 \) we have
\begin{align}
\left( \sum_{n \in \mathbb{N}} \left| \langle A f_{n}, g_{n} \rangle \right|^{2}\right)^{\frac{1}{2}} & = \left( \sum_{n \in \mathbb{N}}  \left| \sum_{k_{1},k_{2} \in \mathbb{Z}^{d}}   \overline{\langle 1, M_{\alpha k_{2}} T_{\alpha k_{1}} \phi \rangle} \langle A_{k_{1},k_{2}} f_{n}, g_{n} \rangle \right|^{2} \right)^{\frac{1}{2}} \nonumber \\
& \leq  \sum_{k_{1},k_{2} \in \mathbb{Z}^{d}}  \left( \sum_{n \in \mathbb{N}} \left| \hat{\phi}(\alpha k_{2}) \right|^{2} \left| \langle A_{k_{1},k_{2}} f_{n}, g_{n} \rangle \right|^{2} \right)^{\frac{1}{2}} \label{eqn:1} \\
& \leq \left( \sum_{k_{2} \in \mathbb{Z}^{d}} \left| \hat{\phi}(\alpha k_{2}) \right| \right) \left( \sup_{k_{2}\in \mathbb{Z}^{d}} \sum_{k_{1}\in \mathbb{Z}^{d}} \left( \sum_{n \in \mathbb{N}}  \left| \langle A_{k_{1},k_{2}} f_{n}, g_{n} \rangle \right|^{2} \right)^{\frac{1}{2}} \right) \nonumber \\
& \leq C \left\| \hat{\phi} \right\|_{W(L^{1}(\mathbb{R}^{d}))}\sup_{k_{2}\in \mathbb{Z}^{d}} \sum_{k_{1}\in \mathbb{Z}^{d}} \left\| A_{k_{1},k_{2}} \right\|_{\mathcal{I}_{2}}, \nonumber
\end{align}
where (\ref{eqn:1}) holds by Minkowski's integral inequality.  Again, by Theorem \ref{thm:intop}, we have
\begin{align*}
\left\|  A_{k_{1},k_{2}} \right\|_{\mathcal{I}_{2}} & \leq \left( \sum_{n_{1},n_{2} \in \mathbb{Z}^{d}}  \sum_{m_{1}, m_{2} \in \mathbb{Z}^{d}} \left| V_{\phi \otimes \phi} k_{k_{1},k_{2}} (\alpha n_{1},\alpha m_{1}, \alpha n_{2},\alpha m_{2}) \right|^{2} \right)^{\frac{1}{2}} \\
& = \left( \sum_{n_{1},n_{2} \in \mathbb{Z}^{d}}  \sum_{m_{1}, m_{2} \in \mathbb{Z}^{d}} \left| V_{\Phi} (b e^{2 \pi i \psi})(\alpha n_{1},\alpha m_{1}, \alpha k_{1}, \alpha n_{2},\alpha m_{2}, \alpha k_{2}) \right|^{2} \right)^{\frac{1}{2}}.
\end{align*}
Thus if
\[ \sup_{k_{2} \in \mathbb{Z}^{d}} \sum_{k_{1} \in \mathbb{Z}^{d}} \left( \sum_{n_{1},n_{2} \in \mathbb{Z}^{d}} \sum_{m_{1}, m_{2}\in \mathbb{Z}^{d}} \left| V_{\Phi} (b e^{2 \pi i \psi})(\alpha n_{1},\alpha m_{1}, \alpha k_{1}, \alpha n_{2},\alpha m_{2}, \alpha k_{2}) \right|^{2} \right)^{\frac{1}{2}} < \infty , \] then  \( A \in \mathcal{I}_{2}(L^{2}(\mathbb{R}^{d})) \).  This quantity is finite if and only if \( b e^{2 \pi i \psi} \in M(c)^{p_{1},p_{2},\dots,p_{6d}} \).

Taking the supremum of \( \sum_{n\in \mathbb{N}} \left| \langle A f_{n}, g_{n} \rangle \right| \) and \( \left( \sum_{n\in \mathbb{N}} \left| \langle A f_{n}, g_{n} \rangle \right|^{2} \right)^{\frac{1}{2}} \) over all orthonormal sequences gives the result for \( p =1 \) and \(p =2 \).  For \( 1 < p < 2 \), the result follows by interpolation.
\end{proof}

\section{Sharp Time-Frequency Conditions on the Symbol of a Fourier Integral Operator with smooth phase}

Theorems \ref{thm:FIOlemma} and \ref{thm:FIOmainone} give conditions on the product of the symbol and phase useful for identifying Schatten class Fourier integral operators. However it would be more convenient to allow the properties of the symbol and properties of the phase separately to determine whether a Fourier integral operator is Schatten class.  Hence, in this section we find conditions on the symbol and phase function of a Fourier integral operator so that their product lies in given mixed modulation spaces.  

We begin by stating some preliminary lemmas.  The following lemma comes from Proposition 3.2 in \cite{boulk}.

\begin{lemma} \label{lemma:boulkdil}  There exists a finite \( C \) such that \[ \left\| \tau \left(t\cdot \right) \right\|_{M^{\infty,1}(\mathbb{R}^{d})} \leq C \left\| \tau \right\|_{M^{\infty,1}(\mathbb{R}^{d})} \hspace{1pc} \forall \tau \in M^{\infty,1}(\mathbb{R}^{d}), t \in [0,1].  \]
\end{lemma}

The next lemma is a special case of Proposition 1.2 in \cite{toftFIO2}, which describes multiplication properties of modulation spaces.

\begin{lemma} \label{lemma:toftmult}   Suppose \( p, q, p_{1}, p_{2}, q_{1}, q_{2} \in [1, \infty] \) satisfy \( \frac{1}{p_{1}} + \frac{1}{p_{2}} = \frac{1}{p} \) and \( \frac{1}{q_{1}} + \frac{1}{q_{2}} = 1 + \frac{1}{q} \).  Then there exists a finite \( C \) such that 
\[ \left\|f g \right\|_{M^{p,q}(\mathbb{R}^{d})} \leq C \left\| f\right\|_{M^{p_{1},q_{1}}(\mathbb{R}^{d})} \left\|g \right\|_{M^{p_{2},q_{2}}(\mathbb{R}^{d})} \hspace{2pc} \forall f \in M^{p_{1},q_{1}}(\mathbb{R}^{d}), g \in M^{p_{2},q_{2}}(\mathbb{R}^{d}). \]  In particular, there is a  finite \( C \) such that
\[ \left\|f g \right\|_{M^{p,q}(\mathbb{R}^{d})} \leq C \left\| f\right\|_{M^{p,q}(\mathbb{R}^{d})} \left\|g \right\|_{M^{\infty,1}(\mathbb{R}^{d})} \hspace{2pc} \forall f \in M^{p,q}(\mathbb{R}^{d}), g \in M^{\infty, 1}(\mathbb{R}^{d}). \] 
\end{lemma}

\begin{corollary} \label{cor:toftmult} There is a finite \( C \) such that
\[ \left\| f^{n} \right\|_{M^{\infty, 1}(\mathbb{R}^{d})} \leq C^{n} \left\| f\right\|^{n}_{M^{\infty,1}(\mathbb{R}^{d})}  \hspace{2pc} \forall f \in M^{\infty,1}(\mathbb{R}^{d}). \] 
\end{corollary}

Although not literally true, functions in \( M^{\infty, 1}(\mathbb{R}^{d}) \) can be thought of as functions in \( L^{\infty}(\mathbb{R}^{d}) \) whose Fourier transform is in \( L^{1}(\mathbb{R}^{d}) \).  Lemma \ref{lemma:toftmult} shows that in some sense multiplying a function \( f \) with an element of \( M^{\infty, 1}(\mathbb{R}^{d}) \) preserves the time-frequency localization of \( f \).    This holds because multiplying \( f \) by an \( L^{\infty}(\mathbb{R}^{d}) \) function preserves the time localization of \( f \), and multiplying \( f \) by a function whose Fourier transform is in \( L^{1}(\mathbb{R}^{d}) \) (i.e., convolving \( \widehat{f} \) with a \( L^{1}(\mathbb{R}^{d}) \) function) preserves the frequency localization of \( f \). This idea extends to the mixed modulation spaces as well, as stated precisely in the following theorem.  The proof is a generalization of the proof of Lemma \ref{lemma:toftmult} (see \cite{toftFIO2} for details), so it is omitted.  

\begin{theorem} \label{thm:generalizeme}  
\begin{itemize}
\item[(a)] Let \( p \in [1,\infty] \).  Suppose \( c \) is a first or second slice permutation and \( p_{1} = p_{2} = \dots = p_{2d} = 2 \), \( p_{2d+1} = p_{2d+2} = \dots = p_{4d} = p \). Then for some finite \( C \) we have for all \(    f \in M(c)^{p_{1},p_{2},\dots, p_{4d} },  g \in M^{\infty,1}(\mathbb{R}^{2d}) \) that \[ \left\| f g \right\|_{M(c)^{p_{1},p_{2},\dots, p_{4d}}} \leq C \left\| f \right\|_{M(c)^{p_{1},p_{2},\dots, p_{4d} }} \left\| g \right\|_{M^{\infty,1}(\mathbb{R}^{2d})} . \]
\item[(b)] Let \( p \in [1,\infty] \).  Suppose \( c \) is a first or second FIO slice permutation and \( p_{1} = p_{2} = \dots = p_{2d} = 2 \), \( p_{2d+1} = p_{2d+2} = \dots = p_{4d} = p \), \( p_{4d+1} = p_{4d+2} = \dots = p_{5d}= 1 \) and \( p_{5d+1} = p_{5d+2} = \dots = p_{6d} = \infty \). Then for some finite \( C \) we have for all \(    f \in M(c)^{p_{1},p_{2},\dots, p_{6d} },  g \in M^{\infty,1}(\mathbb{R}^{3d}) \) that
\[ \left\| f g \right\|_{M(c)^{p_{1},p_{2},\dots, p_{6d}}} \leq C \left\| f \right\|_{M(c)^{p_{1},p_{2},\dots, p_{6d} }} \left\| g \right\|_{M^{\infty,1}(\mathbb{R}^{3d})} . \]
\end{itemize}
\end{theorem}

If \( \varphi \) is real-valued then \( e^{i \varphi} \in L^{\infty} \).  Furthermore \( e^{i \varphi} \) will be smooth only if \( \varphi \) is smooth.  Hence, we expect that \( e^{i \varphi} \in M^{\infty,1} \) if \( \varphi \) is sufficiently smooth.  Theorem \ref{thm:generalizeme} shows that for a Fourier integral operator with smooth phase \( \varphi \), the time-frequency localization of \( ae^{i \varphi} \) is controlled by the time frequency localization of \( a \).  Therefore the conditions in Theorems \ref{thm:FIOlemma} and \ref{thm:FIOmainone} should be satisfied for Fourier integral operators having smooth phase and time-frequency localized symbols.  We make this idea precise in Theorem \ref{thm:compactsupp} and Theorem \ref{thm:othersymbol1}.

\begin{theorem} \label{thm:compactsupp}
\begin{itemize}
\item[(a)] Suppose \(A \) is a Fourier integral operator of the form (\ref{eqn:easyfio}) with symbol \( a \) and phase \( \varphi \).  Let \( c \) be a first or second slice permutation. Fix \( p \in [1,2] \) and let \( p_{1} = p_{2} = \cdots = p_{2d} = 2 \) and \( p_{2d+1} = p_{2d+2} = \cdots = p_{4d} = p \).  If \( a \in M(c)^{p_{1},p_{2},\dots, p_{4d}} \) has compact support and  \( \varphi \in C^{2}\left( \mathbb{R}^{2d} \right) \) is real-valued and satisfies \( D^{\alpha} \varphi \in M^{\infty,1}(\mathbb{R}^{2d}) \) for all multi-indices \( \alpha \) with \( |\alpha| = 2 \), then \( A \in \mathcal{I}_{p}(L^{2}(\mathbb{R}^{d})) \).  

Furthermore, this result is sharp in the sense that if one of the following conditions holds, then there are Fourier integral operators not in \( \mathcal{I}_{p}(L^{2}(\mathbb{R}^{d})) \) with symbols in \( M(c)^{q_{1},q_{2},\dots, q_{4d}} \)  and phase functions \(\varphi \) satisfying \( \varphi \in C^{2}(\mathbb{R}^{2d}) \) and  \( D^{\alpha} \varphi \in M^{\infty,1}(\mathbb{R}^{2d}) \) for all multi-indices  \( \alpha \) with  \( |\alpha| = 2 \).
\begin{itemize}
\item[(i)]  At least one of \(  q_{1} ,  \dots  , q_{2d} \) is larger than 2.
\item[(ii)] At least one of \(  q_{2d+1} , \dots , q_{4d} \) is larger than \( p \).
\end{itemize} 

\item[(b)] Suppose \(A \) is a Fourier integral operator of the form (\ref{eqn:hardfio}) with symbol \( b \) and phase \( \psi \).  Let \( c \) be a first or second FIO slice permutation. Fix \( p \in [1,2] \) and let \( p_{1} = p_{2} = \cdots = p_{2d} = 2 \), \( p_{2d+1} = p_{2d+2} = \cdots = p_{4d} = p \), \( p_{4d+1} = p_{4d+2} = \dots = p_{5d}= 1 \) and \( p_{5d+1} = p_{5d+2} = \dots = p_{6d} = \infty \).  If \( b \in M(c)^{p_{1},p_{2},\dots, p_{6d}} \) has compact support and  \( \psi \in C^{2}\left( \mathbb{R}^{3d} \right) \) is real valued and satisfies \( D^{\alpha} \psi \in M^{\infty,1}(\mathbb{R}^{3d}) \) for all multi-indices \( \alpha \) with \( |\alpha| = 2 \), then \( A \in \mathcal{I}_{p}(L^{2}(\mathbb{R}^{d})) \).

 Furthermore, this result is sharp in the sense that if one of the following conditions holds, then there are Fourier integral operators that are not in   \( \mathcal{I}_{p}(L^{2}(\mathbb{R}^{d})) \) with symbols in \(  M(c)^{q_{1},q_{2},\dots, q_{5d}, p_{5d+1}, p_{5d+2}, \dots p_{6d}} \) and phase functions \(\psi \) satisfying \( \psi \in C^{2}(\mathbb{R}^{3d}) \) and  \( D^{\alpha} \psi \in M^{\infty,1}(\mathbb{R}^{3d}) \) for all multi-indices  \( \alpha \) with  \( |\alpha| = 2 \).
\begin{itemize}
\item[(i)] At least one of \( q_{1}, q_{2},\dots, q_{2d} \) is larger than 2.
\item[(ii)]  At least one of \( q_{2d+1}, q_{2d+2}, \dots, q_{4d} \) is larger than \( p \).
\item[(iii)] At least one of \( q_{4d+1},q_{4d+2},\dots,q_{5d} \) is larger than 1.
\end{itemize}

 \end{itemize}
\end{theorem}  
\begin{proof} We prove (b).  Statement (a) can be proven similarly.

By Theorem \ref{thm:FIOmainone}, it suffices to prove \( \left\| be^{2 \pi i \psi} \right\|_{M(c)^{p_{1},\dots, p_{6d}}} < \infty \). Write \( \psi = \psi_{1} + \psi_{2} \), where
\[ \psi_{1}(w) =  \varphi(0,0,0) + \sum_{ \left| \alpha \right| = 1 } \left( D^{\alpha} \varphi \right)(0,0,0) w^{\alpha} \] and 
\[ \psi_{2}(w) = \sum_{ \left| \alpha \right| = 2 } \frac{2}{\alpha!} \left( \int_{0}^{1} \left( 1-t \right) \left( D^{\alpha} \varphi \right)(tw) \, \ud t \right) w^{\alpha}.\]
Choose \( \CHI \in C_{c}^{\infty}(\mathbb{R}^{3d}) \)  such that \( \CHI(w) = 1  \) for all \( w \) in the support of \( b \).  Then
\begin{align*}
\left\| be^{2 \pi i \psi} \right\|_{M(c)^{p_{1},\dots, p_{6d}}} & = \left\| be^{2 \pi i \psi_{1}} e^{2 \pi i \psi_{2}} \right\|_{M(c)^{p_{1},\dots, p_{6d}}} \nonumber \\
& = \left\| be^{2 \pi i \psi_{1}} e^{2 \pi i \CHI \psi_{2}}\right\|_{M(c)^{p_{1},\dots, p_{6d}}} \nonumber \\
& \leq \left\| be^{2 \pi i \psi_{1}} \right\|_{M(c)^{p_{1},\dots, p_{6d}}} \left \| e^{2 \pi i \CHI \psi_{2}}\right\|_{M^{\infty,1}} \\
& < \infty.
\end{align*}
 Notice that
\[ be^{2 \pi i \psi_{1}} = M_{q} \left(e^{2 \pi i \varphi(0,0,0)} b \right), \]
where the components of \( q \in \mathbb{R}^{3d} \) are \(  \left( D^{\alpha} \varphi \right)(0,0,0) \) for multi-indices \( \alpha \) with  \( \left| \alpha \right| = 1  \).  Since \( b \in M(c)^{p_{1},p_{2},\dots,p_{6d}} \), we have  \( \left\| be^{2 \pi i \psi_{1}} \right\|_{M(c)^{p_{1},\dots, p_{6d}}} = \left\| a\right\|_{M(c)^{p_{1},\dots, p_{6d}}} < \infty \) as well.  

Now we show that \( e^{2 \pi i \CHI \psi_{2}} \in M^{\infty,1}(\mathbb{R}^{3d}) \).  To this end, choose finite \( C \) such that \[ \left\|f g \right\|_{M^{\infty,1}(\mathbb{R}^{d})} \leq C \left\| f\right\|_{M^{\infty, 1}(\mathbb{R}^{d})} \left\|g \right\|_{M^{\infty,1}(\mathbb{R}^{d})} \hspace{2pc} \forall f, g \in M^{\infty, 1}(\mathbb{R}^{d}). \] We have
\begin{align*}
\left\| e^{2 \pi i \CHI \psi_{2}} \right\|_{M^{\infty,1}} & = \left\|\sum_{n \geq 0} \frac{\left( 2 \pi i \CHI \psi_{2}\right)^{n}}{n!} \right\|_{M^{\infty,1}} \nonumber \\
& \leq \sum_{n \geq 0} \frac{\left\| \left(2 \pi i \CHI \psi_{2}\right)^{n} \right\|_{M^{\infty,1}}}{n!}  \nonumber \\
& \leq \sum_{n \geq 0} \frac{\left( 2 \pi  C\right)^{n} \left\| \CHI \psi_{2} \right\|^{n}_{M^{\infty,1}}}{n!} \label{eqn:compactsymb2} \\
& = e^{2 \pi C \left\| \CHI \psi_{2} \right\|_{M^{\infty,1}}}. \nonumber
\end{align*}

By Lemma \ref{lemma:boulkdil}, we can choose \( C' \) so that \[ \left\| \tau \left(t\cdot \right) \right\|_{M^{\infty,1}(\mathbb{R}^{d})} \leq C' \left\| \tau \right\|_{M^{\infty,1}(\mathbb{R}^{d})} \hspace{1pc} \forall \tau \in M^{\infty,1}(\mathbb{R}^{d}), t \in [0,1].  \]  Since \( \CHI \in C_{c}^{\infty}(\mathbb{R}^{d}) \), we have
\( \CHI(w)  w^{\alpha} \in M^{\infty,1}(\mathbb{R}^{d}) \) for all \( \alpha \) with \(|\alpha| = 2. \)  Thus
\begin{align*}
& \left\| \CHI \psi_{2} \right\|_{M^{\infty,1}} \nonumber \\
 & =  \left\| \sum_{ \left| \alpha \right| = 2 } \CHI(w) w^{\alpha} \frac{2}{\alpha!} \left( \int_{0}^{1} \left( 1-t \right) \left( D^{\alpha} \varphi \right)(tw) \, \ud t \right) \right\|_{M^{\infty,1}}  \nonumber \\
& \leq   \sum_{ \left| \alpha \right| = 2 } \frac{2}{\alpha!}  \left\|\CHI(w) w^{\alpha} \left( \int_{0}^{1} \left( 1-t \right) \left( D^{\alpha} \varphi \right)(tw) \, \ud t \right) \right\|_{M^{\infty,1}}  \nonumber \\
& \leq   \sum_{ \left| \alpha \right| = 2 } \frac{2C}{\alpha!}  \left\|\CHI(w) w^{\alpha} \right\|_{M^{\infty,1}}   \left\| \int_{0}^{1} \left( 1-t \right) \left( D^{\alpha} \varphi \right)(tw) \, \ud t  \right\|_{M^{\infty,1}}  \label{eqn:compactsymb3} \\
& = \sum_{ \left| \alpha \right| = 2 } \frac{2C}{\alpha!}  \left\|\CHI(w) w^{\alpha} \right\|_{M^{\infty,1}}   \int \sup_{x \in \mathbb{R}^{d}} \left| \int \int_{0}^{1} \left( 1-t \right) \left( D^{\alpha} \varphi \right)(tw) \overline{ M_{\xi} T_{x} \phi(w) } \, \ud t \, \ud w \right| \, \ud \xi \nonumber \\
& \leq \sum_{ \left| \alpha \right| = 2 } \frac{2C}{\alpha!}  \left\|\CHI(w) w^{\alpha} \right\|_{M^{\infty,1}}   \int \sup_{x \in \mathbb{R}^{d}}  \int_{0}^{1} \left( 1-t \right) \left| \int \left( D^{\alpha} \varphi \right)(tw) \overline{ M_{\xi} T_{x} \phi(w) }  \, \ud w \right| \, \ud t \, \ud \xi \nonumber \\
& \leq \sum_{ \left| \alpha \right| = 2 } \frac{2C}{\alpha!}  \left\|\CHI(w) w^{\alpha} \right\|_{M^{\infty,1}} \int_{0}^{1} \left( 1-t \right)  \int \sup_{x \in \mathbb{R}^{d}}   \left| \int \left( D^{\alpha} \varphi \right)(tw) \overline{ M_{\xi} T_{x} \phi(w) }  \, \ud w \right| \, \ud \xi  \, \ud t \nonumber \\
& =   \sum_{ \left| \alpha \right| = 2 } \frac{2C}{\alpha!}  \left\|\CHI(w) w^{\alpha} \right\|_{M^{\infty,1}}    \int_{0}^{1} \left( 1-t \right) \left\| \left(D^{\alpha} \varphi \right)(tw) \right\|_{M^{\infty,1}} \, \ud t   \nonumber   \\
& \leq   \sum_{ \left| \alpha \right| = 2 } \frac{2C}{\alpha!}  \left\|\CHI(w) w^{\alpha} \right\|_{M^{\infty,1}}    \int_{0}^{1} \left( 1-t \right) C' \left\| \left(D^{\alpha} \varphi \right)(w) \right\|_{M^{\infty,1}} \, \ud t    \nonumber \\
& \leq   \sum_{ \left| \alpha \right| = 2 } \frac{2CC'}{\alpha!}  \left\|\CHI(w) w^{\alpha} \right\|_{M^{\infty,1}}  \left\| \left(D^{\alpha} \varphi \right)(w) \right\|_{M^{\infty,1}}    \nonumber \\
& < \infty \nonumber. 
\end{align*}

Now we will prove sharpness.  Let \( c_{1} \) be the permutation of \( \left\{ 1,2, \dots, 4d \right\} \) such that 
\[ c_{1}(1)  = c(1), c_{1}(2)  = c(2), \dots , c_{1}(d)  = c(d,) \]
\[ c_{1}(d+1)  = c(d+1), c_{1}(d+2)  = c(d+2), \dots , c_{1}(2d)  = c(2d), \]
\[ c_{1}(2d+1)  = c(3d+1), c_{1}(2d+2)  = c(3d+2),  \dots , c_{1}(3d) = c(4d), \] and
\[ c_{1}(3d+1)  = c(4d+1), c_{1}(3d+2) = c(4d+2), \dots, c_{1}(4d)  = c(5d), \]
and let \( c_{2} \) be the permutation of \( \left\{ 1,2,\dots, 2d \right\} \) such that
\[c_{2}(1)  = c(2d+1) -4d, c_{2}(2) = c(2d+2) -4d, \dots , c_{2}(d)   = c(3d) - 4d, \] and
\[ c_{2}(d+1)  = c(5d+1) - 4d, c_{2}(d+2)  = c(5d+2)-4d,  \dots, c_{2}(2d)  = c(6d) - 4d. \]
Notice that \( c_{1} \) is a first or second slice permutation.

Let \( b(x,y,\xi) = b_{1}(x,y) b_{2}(\xi) \).  Then
\[ \left\| b \right\|_{M(c)^{q_{1},q_{2},\dots, q_{5d}, p_{5d+1}, p_{5d+2}, \dots p_{6d}} } = \left\| b_{1} \right\|_{M(c_{1})^{q_{1}, q_{2}, \dots, q_{4d}}} \left\| b_{2} \right\|_{ M(c_{2})^{q_{4d+1}, \cdots, q_{5d}, p_{5d+1}, \dots, p_{6d}}}. \]  Let \( A \) be the Fourier integral operator of the form (\ref{eqn:hardfio})  with symbol \( b \) and phase function \( \psi = 1 \).  Then
\[ Af(x) = \int b_{2}(\xi) \, \ud \xi \, \int b_{1}(x,y) f(y) \, \ud y. \]  If (iii) holds, then choose  \( b_{1} \in M(c_{1})^{q_{1}, q_{2}, \dots, q_{4d}} \) and \( b_{2} \in M(c_{2})^{q_{4d+1}, \cdots, q_{5d}, p_{5d+1}, \dots, p_{6d}}  \) such that \( \int b_{2}(\xi) \, \ud \xi = \infty.\)  Then \( A:L^{2}(\mathbb{R}^{d}) \to L^{2}(\mathbb{R}^{d}) \) is not well-defined and hence not in \( \mathcal{I}_{p}(L^{2}(\mathbb{R}^{d})) \).  Otherwise, if (i) or (ii) hold, then choose \( b_{1} \in M(c_{1})^{q_{1}, q_{2}, \dots, q_{4d}} \) and \( b_{2} \in M(c_{2})^{q_{4d+1}, \cdots, q_{5d}, p_{5d+1}, \dots, p_{6d}}  \) such that \(  \int b_{2}(\xi) \, \ud \xi  \neq 0 \) and the integral operator with kernel \( b_{1} \) is not in \( \mathcal{I}_{p}(L^{2}(\mathbb{R}^{d})) \) (such a choice is possible by Theorem \ref{thm:intop}).  Then \( A \notin \mathcal{I}_{p}(L^{2}(\mathbb{R}^{d})) \).
\end{proof}

In the remainder of this section, we develop alternative conditions on the symbol and phase function of a Fourier integral operator so that their product lies in mixed modulation spaces  relevant to Schatten class integral operators.   First, a technical lemma is needed.  See ??? \cite{grochbook} for a proof.

\begin{lemma} \label{lemma:536} Suppose \( \Phi \in M^{1,1}\left(\mathbb{R}^{d} \right) \) and \( M \) is a $d \times d$ self-adjoint matrix.  Define a operator \( S_{M} \) by
\[  S_{M} f(w) = e^{\pi i w \cdot M w } f(w), \hspace{2pc} \forall f \in M^{\infty, \infty}\left(\mathbb{R}^{d} \right) .\]  Then 
\[ \left| V_{\Phi} S_{M} f \left(x, \xi \right)   \right| = \left| V_{S_{-M}\Phi}  f \left(x, \xi - Mx \right)   \right|, \hspace{2pc} \forall x,\xi \in \mathbb{R}^{d}.  \]
\end{lemma}

\begin{theorem} \label{thm:othersymbol1}  Let \( p \in [1, 2] \).  
\begin{itemize}
\item[(a)] Suppose \(A \) is a Fourier integral operator of the form (\ref{eqn:easyfio}) with symbol \( a \) and phase \( \varphi \).  Let \( c \) be a first or second slice permutation and let \( p_{1} = p_{2} = \cdots = p_{2d} = 2 \) and \( p_{2d+1} = p_{2d+2} = \cdots = p_{4d} = p \).  Suppose \( a \in M(c)^{p_{1},p_{2},\dots, p_{4d}} \) and  \( \varphi \in C^{2}\left( \mathbb{R}^{2d} \right) \) has constant second order partial derivatives with \( \varphi_{x_{i}y_{j}} = 0 \) for all \( i,j \in \left\{ 1,2, \dots , d \right\} \).  Then \( A \in \mathcal{I}_{p}(L^{2}(\mathbb{R}^{d})) \).  

Furthermore, this result is sharp in the sense that if one of the following conditions hold then there are Fourier integral operators not in \( \mathcal{I}_{p}(L^{2}(\mathbb{R}^{d})) \) with symbols in \( M(c)^{q_{1},q_{2},\dots, q_{4d}} \)  and phase functions \( \varphi \in C^{2}(\mathbb{R}^{2d}) \) with constant second order partial derivatives.
\begin{itemize}
\item[(i)]  At least one of \(  q_{1} ,  \dots  , q_{2d} \) is larger than 2.
\item[(ii)] At least one of \(  q_{2d+1} , \dots , q_{4d} \) is larger than \( p \).
\end{itemize} 
\item[(b)] Suppose \(A \) is a Fourier integral operator of the form (\ref{eqn:hardfio}) with symbol \( b \) and phase \( \psi \).  Let  \( c \) be a first or second FIO symbol permutation and \( p_{1} = \dots = p_{d} = \infty \), \( p_{d+1} = p_{d+ 2} = \dots = p_{3d} = 2 \), \( p_{3d+1} = p_{3d+2} = \dots = p_{5d} = p \) and \( p_{5d+1} = p_{5d+2} = \dots p_{6d}= 1 \). Suppose \( b \in M(c)^{p_{1}, p_{2}, \dots, p_{6d}} \), \( \psi \in C^{2}(\mathbb{R}^{3d}) \), all the second order partial derivatives of \( \psi \) are constant and \( \psi_{x_{i}y_{j}} = 0 \) for all \( i,j \in \left\{ 1,2, \dots , d \right\} \). Then \( A \in \mathcal{I}_{p}(L^{2}(\mathbb{R}^{d})) \).  

Furthermore, this result is sharp in the sense that if one of the following conditions holds, then there are Fourier integral operators that are not in   \( \mathcal{I}_{p}(L^{2}(\mathbb{R}^{d})) \) with symbols in \(  M(c)^{p_{1},\cdots,p_{d}, q_{d+1}, \cdots q_{6d}} \) and phase functions  \( \psi \in C^{2}(\mathbb{R}^{3d}) \)  with constant second order partial derivatives.
\begin{itemize}
\item[(i)] At least one of \( q_{d+1}, q_{d+2},\dots, q_{3d} \) is larger than 2.
\item[(ii)]  At least one of \( q_{3d+1}, q_{3d+2}, \dots, q_{5d} \) is larger than \( p \).
\item[(iii)] At least one of \( q_{5d+1},q_{5d+2},\dots,q_{6d} \) is larger than 1.
\end{itemize}
\end{itemize}
\end{theorem}
\begin{proof}  We prove (b) in the case \( c \) is a first FIO symbol permutation.  The other parts and cases can be proved similarly.

Let \( c' \) be a first FIO slice permutation.  Set \( r_{1} = r_{2} = \cdots = r_{2d} = 2 \), \( r_{2d+1} = r_{2d+2} = \cdots = r_{4d} = p \), \( r_{4d+1} = r_{4d+2} = \dots = r_{5d}= 1 \) and \( r_{5d+1} = r_{5d+2} = \dots = r_{6d} = \infty \).  By Theorem \ref{thm:FIOmainone}, it suffices to show \( b e^{2 \pi i \psi} \in M(c')^{r_{1},r_{2},\dots, r_{6d}} \).

Write \( \psi = \psi_{1} + \psi_{2} \), where
\[ \psi_{1}(w) =  \varphi(0,0,0) + \sum_{ \left| \alpha \right| = 1 } \left( D^{\alpha} \varphi \right)(0,0,0) w^{\alpha} \] and 
\[ \psi_{2}(w) = \sum_{ \left| \alpha \right| = 2 } \frac{2}{\alpha!} \left( \int_{0}^{1} \left( 1-t \right) \left( D^{\alpha} \varphi \right)(tw) \, \ud t \right) w^{\alpha}.\]
Notice that  \( e^{2 \pi i \psi_{2}(w)} = e^{i \pi w \cdot M w} \) where \( M \) is the block matrix
\[ M = 2 \pi \left[ \begin{array}{ccc}   M_{1} & M_{2} & M_{3} \\
  M^{*}_{2} & M_{4} & M_{5} \\
   M^{*}_{3} & M^{*}_{5} & M_{6}
  \end{array}  \right],\]   with 
  \[ \left( M_{1} \right)_{i,j} = \frac{\psi_{x_{i}x_{j}}(0,0,0)}{2}  \hspace{2pc} \forall i,j \in \left\{ 1,2, \dots , d \right\}, \]
    \[ \left( M_{2} \right)_{i,j} = \frac{\psi_{x_{i}y_{j}}(0,0,0)}{2} = 0 \hspace{2pc} \forall i,j \in \left\{ 1,2, \dots , d \right\}, \]
      \[ \left( M_{3} \right)_{i,j} = \frac{\psi_{x_{i}\xi_{j}}(0,0,0)}{2}  \hspace{2pc} \forall i,j \in \left\{ 1,2, \dots , d \right\}, \]
        \[ \left( M_{4} \right)_{i,j} = \frac{\psi_{y_{i}y_{j}}(0,0,0)}{2}  \hspace{2pc} \forall i,j \in \left\{ 1,2, \dots , d \right\}, \]
          \[ \left( M_{5} \right)_{i,j} = \frac{\psi_{y_{i}\xi_{j}}(0,0,0)}{2}  \hspace{2pc} \forall i,j \in \left\{ 1,2, \dots , d \right\}, \] and
            \[ \left( M_{6} \right)_{i,j} = \frac{\psi_{\xi_{i}\xi_{j}}(0,0,0)}{2}  \hspace{2pc} \forall i,j \in \left\{ 1,2, \dots , d \right\}. \]

            Thus
           \begin{align} & \left|V_{\Phi}\left( be^{2 \pi i \psi} \right)(x_{1}, x_{2}, x_{3},\xi_{1}, \xi_{2}, \xi_{3}) \right| \nonumber \\
            & \indent \indent = \left| V_{\Phi}\left( be^{2 \pi i \psi} \right)(x,\xi) \right| \nonumber \\
           &  \indent \indent = \left|V_{\Phi}S_M \left(be^{2 \pi i \psi_{1}}\right) (x,\xi) \right| \nonumber \\
             & \indent \indent = \left|V_{S_{-M} \Phi}\left(be^{2 \pi i \psi_{1}}\right) (x,\xi - Mx) \right| \label{eqn:cite536} \\
               & \indent \indent = \left|V_{S_{-M} \Phi}\left(be^{2 \pi i \psi_{1}}\right)(x_{1}, x_{2}, x_{3}, \xi_{1} - M_{1} x_{1} - M_{2} x_{2} - M_{3} x_{3}, \right. \nonumber \nonumber \\
          & \indent \indent \indent \indent       \left. \xi_{2} - M^{*}_{2} x_{1} - M_{4} x_{2} -M_{5} x_{3}, \xi_{3} - M^{*}_{3} x_{1} - M^{*}_{5} x_{2} - M_{6} x_{3}) \right| \nonumber \\
          & \indent \indent = \left|V_{S_{-M} \Phi}\left(b e^{2 \pi i \psi_{1}}\right)(x_{1}, x_{2}, x_{3}, \xi_{1} - M_{1} x_{1} - M_{3} x_{3}, \right.\nonumber \\
          & \indent \indent \indent \indent       \left. \xi_{2} - M_{4} x_{2} -M_{5} x_{3}, \xi_{3} - M^{*}_{3} x_{1} - M^{*}_{5} x_{2} - M_{6} x_{3}) \right|  \nonumber
           \end{align}
where (\ref{eqn:cite536}) follows from Lemma \ref{lemma:536}.

          Hence
           \begin{align*}
           & \left\|be^{2 \pi i \psi} \right\|_{M(c')^{r_{1}, \dots, r_{6d}}} \\
           & = \sup_{\xi_{3}} \int \left( \iint \left( \iint \left| V_{\Phi}( be^{2 \pi i \psi}  ) ( x_{1}, x_{2}, x_{3}, \xi_{1}, \xi_{2}, \xi_{3}) \right|^{2}  \, \ud \xi_{1} \, \ud x_{1} \right)^{\frac{p}{2}}  \, \ud \xi_{2} \, \ud x_{2}  \right)^{\frac{1}{p}} \, \ud x_{3} \\
            & = \sup_{\xi_{3}} \int \left( \iint \left( \iint  \left|V_{S_{-M} \Phi}\left(be^{2 \pi i \psi_{1}}\right)(x_{1}, x_{2}, x_{3}, \xi_{1} - M_{1} x_{1} - M_{3} x_{3},  \right. \right. \right. \\
            & \indent \indent  \left. \left. \left. \xi_{2} - M_{4} x_{2} -M_{5} x_{3}, \xi_{3} - M^{*}_{3} x_{1} - M^{*}_{5} x_{2} - M_{6} x_{3}) \right|^{2}  \, \ud \xi_{1} \, \ud x_{1} \right)^{\frac{p}{2}}   \, \ud \xi_{2} \, \ud x_{2} \right)^{\frac{1}{p}} \, \ud x_{3} \\
            & \leq  \int \left( \iint \left( \iint \sup_{\xi_{3}} \left|V_{S_{-M} \Phi}\left(be^{2 \pi i \psi_{1}}\right)(x_{1}, x_{2}, x_{3}, \xi_{1} - M_{1} x_{1} - M_{3} x_{3},  \right. \right. \right. \\
            & \indent \indent    \left. \left. \left.  \xi_{2} - M_{4} x_{2} -M_{5} x_{3}, \xi_{3} - M^{*}_{3} x_{1} - M^{*}_{5} x_{2} - M_{6} x_{3}) \right|^{2}  \, \ud \xi_{1} \, \ud x_{1} \right)^{\frac{p}{2}}   \, \ud \xi_{2} \, \ud x_{2} \right)^{\frac{1}{p}} \, \ud x_{3} \\
            & =  \int \left( \iint \left( \iint \sup_{\xi_{3}} \left| V_{S_{-M}\Phi}( be^{2 \pi i \psi_{1}}  ) ( x_{1}, x_{2}, x_{3}, \xi_{1}, \xi_{2}, \xi_{3}) \right|^{2}  \, \ud \xi_{1} \, \ud x_{1} \right)^{\frac{p}{2}}  \, \ud \xi_{2} \, \ud x_{2}  \right)^{\frac{1}{p}} \, \ud x_{3} \\
            & \equiv \left\| be^{2 \pi i \psi_{1}} \right\|_{M(c)^{p_{1}, \dots, p_{6d}}}.
           \end{align*}
           
           As in the proof of Theorem \ref{thm:compactsupp}, we have \[ be^{2 \pi i \psi_{1}} = M_{q} \left(e^{2 \pi i \varphi(0,0,0)} b \right), \]
where the components of \( q \in \mathbb{R}^{3d} \) are \(  \left( D^{\alpha} \varphi \right)(0,0,0) \) for multi-indices \( \alpha \) with  \( \left| \alpha \right| = 1  \).  Therefore   
           \[  \left\| be^{2 \pi i \psi_{1}} \right\|_{M(c)^{p_{1}, \dots, p_{6d}}} =  \left\| b \right\|_{M(c)^{p_{1}, \dots, p_{6d}}} < \infty, \]
   which implies        \( b e^{2 \pi i \psi} \in  M(c')^{r_{1},r_{2},\dots,r_{6d}} \), as desired.
   
   Sharpness can be proved similarly to the proof of sharpness in Theorem \ref{thm:compactsupp}.
\end{proof}

In the literature on Fourier integral operators, it is common to avoid degenerate operators by assuming
\[ \left| \textrm{det} \left( \varphi_{xy} \right) \right| \geq d > 0 \] for Fourier integral operators of the form (\ref{eqn:easyfio}) and
\[ \left| \textrm{det} \left( \begin{array}{cc} 
\psi_{xy} & \psi_{x \xi} \\ 
 \psi_{y \xi} & \psi_{\xi \xi}
\end{array} \right) \right| \geq d > 0 \]
for Fourier integral operators of the form (\ref{eqn:hardfio}).  Although Theorem \ref{thm:othersymbol1} is only applicable to degenerate Fourier integral operators of the form (\ref{eqn:easyfio}) and certain types of both degenerate and non-degenerate  Fourier integral operators of the form (\ref{eqn:hardfio}), we can generalize the proof technique to apply to both degenerate and non-generate Fourier integral operators of type  (\ref{eqn:easyfio}) and  (\ref{eqn:hardfio}).

\begin{theorem} \label{thm:nodegen}  Let \( p \in [1, 2] \).  
\begin{itemize}
\item[(a)] Suppose \(A \) is a Fourier integral operator of the form (\ref{eqn:easyfio}) with symbol \( a \) and phase \( \varphi \).  Let \( c \) be a permutation on \( \left\{ 1,2,\cdots,4d \right\} \) satisfying one of the following conditions. 
\begin{itemize}
\item[(i)] \( c \) maps \( \left\{3d+1,\cdots,4d \right\} \) to \( \left\{ 1,2,\cdots,d \right\} \) or
\item[(ii)] \( c \) maps \( \left\{2d+1,\cdots,3d \right\} \) to \( \left\{ 1,2,\cdots,d \right\} \)
\end{itemize}  
Let \( p_{1} = p_{2} = \cdots = p_{d} = 2 \) and \( p_{d+1} = p_{d+2} = \cdots = p_{4d} = p \).  Suppose \( a \in M(c)^{p_{1},p_{2},\dots, p_{4d}} \) and  \( \varphi \in C^{2}\left( \mathbb{R}^{2d} \right) \) has constant second order partial derivatives. Then \( A \in \mathcal{I}_{p}(L^{2}(\mathbb{R}^{d})) \).  
 
\item[(b)] Suppose \(A \) is a Fourier integral operator of the form (\ref{eqn:hardfio}) with symbol \( b \) and phase \( \psi \).   Let \( c \) be a permutation on \( \left\{ 1,2,\cdots,6d \right\} \) that maps \( \left\{5d+1,\cdots,6d \right\} \) to \( \left\{ 1,2,\cdots,d \right\} \),  maps \( \left\{2d+1,\cdots,3d \right\} \) to \( \left\{ 5d+1,\cdots,6d \right\} \), and  satisfies one of the following conditions. 
\begin{itemize}
\item[(i)] \( c \) maps \( \left\{4d+1,\cdots,5d \right\} \) to \( \left\{ d+1,d+2,\cdots,2d \right\} \) or
\item[(ii)] \( c \) maps \( \left\{3d+1,\cdots,4d \right\} \) to \( \left\{ d+1,d+2,\cdots,2d \right\} \)
\end{itemize}   Let  \( p_{1} = \dots = p_{d} = \infty \), \( p_{d+1} = p_{d+ 2} = \dots = p_{2d} = 2 \), \( p_{2d+1} = p_{3d+2} = \dots = p_{5d} = p \) and \( p_{5d+1} = p_{5d+2} = \dots p_{6d}= 1 \). Suppose \( b \in M(c)^{p_{1}, p_{2}, \dots, p_{6d}} \), \( \psi \in C^{2}(\mathbb{R}^{3d}) \), all the second order partial derivatives of \( \psi \) are constant. Then \( A \in \mathcal{I}_{p}(L^{2}(\mathbb{R}^{d})) \).  
\end{itemize}
\end{theorem}
\begin{proof}
We prove (a) in the case that (i) holds.  The other case and statement (b) are proved similarly. 

Notice that by By Lemma 4.7  in \cite{mine}, for any permutation \( c' \) on \( \left\{ 1,\cdots, 4d \right\} \) and any exponents \( q_{i}, r_{i} \in [1, \infty] \) satisfying \( q_{i} \leq r_{i} \) for all \( i \in  \left\{ 1,\cdots, 4d \right\} \) we have \( M(c')^{q_{1}, \cdots, q_{4d}}  \subset  M(c')^{r_{1}, \cdots, r_{4d}} \).

 Choose a permutation \( c' \) on \( \left\{ 1,\cdots, 4d \right\} \) such that  \( c' \) maps  \( \left\{3d+1,\cdots,4d \right\} \) to \( \left\{ 1,2,\cdots,d \right\} \),  \( \left\{d+1,\cdots,2d \right\} \) to \( \left\{ d+1,\cdots,2d \right\} \), and  \( \left\{1, \cdots, d, 2d+1, \dots, 3d \right\} \) to \( \left\{ 2d+1,\cdots,4d \right\} \). Notice that \( c' \) is a second slice permutation.

By Theorem \ref{thm:FIOlemma}, it suffices to show \( a e^{2 \pi i \varphi} \in M(c')^{q_{1}, \cdots, q_{4d} } \), where \( q_{1} = \dots = q_{2d} = 2 \) and  where \( q_{2d+1} = \dots = q_{4d} = p \). Using Lemma \ref{lemma:536} as in the proof of Theorem \ref{thm:othersymbol1}, we can show \( a e^{2 \pi i \varphi} \in M(c)^{p_{1}, \cdots, p_{4d}} \), and by our choice of \( c' \) we have \( M(c)^{p_{1}, \cdots, p_{4d}} =  M(c')^{p_{1}, \cdots, p_{4d}} \).  By Lemma 4.7 in \cite{mine}, \(  M(c')^{p_{1}, \cdots, p_{4d}} \subset M(c')^{q_{1}, \cdots, q_{4d} } \) since \( p_{d+1} \leq  q_{d+1}, \dots, p_{2d} \leq q_{2d} \).   Thus \( a e^{2 \pi i \varphi} \in M(c')^{q_{1}, \cdots, q_{4d} } \).
\end{proof}




\section*{Acknowledgements}
The author was partially supported by NSF Grant DMS-0806532.

\bibliographystyle{plain}
\bibliography{FIOpubD4}


\end{document}